\documentclass[10pt]{article}
\usepackage[latin1]{inputenc}
\usepackage{graphics}
\usepackage{amssymb}
\usepackage{amsfonts}
\usepackage[thmmarks]{ntheorem}
\usepackage[all]{xy}

\newcommand{\eqref}[1]{(\ref{#1})}

\def\RR{{\mathbb R}}
\def\CC{{\mathbb C}}
\def\NN{{\mathbb N}}

\def\dsf{\mbox{\rm D.S.F}}

\def\End{\mbox{\rm End}}

\def\bra{\langle}
\def\ket{\rangle}

\def\bea{\begin{eqnarray}}
\def\eea{\end{eqnarray}}

\def\be{\begin{equation}}
\def\ee{\end{equation}}

\newtheorem{theorem}{Theorem}
\newtheorem{definition}{Definition}
\newtheorem{lemma}{Lemma}

\newtheorem{propo}{Proposition}
\theoremstyle{nonumberplain}
\theorembodyfont{\normalfont}
\theoremseparator{:}
\theoremsymbol{$\P$}

\newtheorem{demo}{Proof}
\begin{document}
\title{Doppler shift in semi-Riemannian signature and the non-uniqueness of the Krein space of spinors}
\author{Fabien Besnard\footnote{P\^ole de recherche M.L. Paris, EPF, 3 bis rue Lakanal, F-92330 Sceaux.}, Nadir Bizi\footnote{Institut de Min\'eralogie, de Physique des Mat\'eriaux et de
Cosmochimie,  Sorbonne Universit\'es, UMR CNRS 7590, 
Universit\'e Pierre et Marie Curie-Paris 06,
  Mus\'eum National d'Histoire Naturelle, IRD UMR 206,
4 place Jussieu, F-75005 Paris, France}\\
{\small fabien.besnard@epf.fr, nadir.bizi@impmc.upmc.fr}}
\maketitle
\begin{abstract} We give     examples illustrating the fact that the different space/time splittings of the tangent bundle of a semi-Riemannian spin manifold give rise to non-equivalent norms on the space of compactly supported sections of the spinor bundle, and as a result, to different completions. We give a necessary and sufficient condition for two space/time splittings to define equivalent norms in terms of a generalized Doppler shift between maximal negative definite   subspaces. We explore some consequences for the Noncommutative Geometry program.
\end{abstract}

\section{Introduction}
Since the seminal work   \cite{baum}, it is well-known that the space $\Gamma_c(S)$ of compactly supported sections of the spinor bundle, or spinor fields, of a semi-riemannian spin manifold $(M,g)$ can be completed into a Krein space. This is done thanks to a splitting of the tangent bundle of $M$ into spacelike and timelike subbundles. More precisely, there exists on $\Gamma_c(S)$ a hermitian form $(.,.)$ of indefinite signature which is uniquely defined up to a local rescaling, and the splitting can be used to turn it into a scalar product. This scalar product can then be used to complete $\Gamma_c(S)$ into a Krein space. In the case where $M$ is non-compact and $g$ is non-Riemannian, the completion will generally depend on the splitting. Yet this dependence has not been explored or even emphasized in the literature (to the authors knowledge). There has even been at least one opposite claim (sadly, by authors themselves, in \cite{part1}, p. 17). This question is however of uttermost importance for the formulation of semi-Riemannian Noncommutative Geometry. Indeed, in the mainstream approach to this subject ( \cite{stro}, \cite{PS}, \cite{vddgpr},  \cite{francoeckstein}, \cite{SST2}) one defines a semi-Riemannian replacement for spectral triples, and though the precise axioms differ, they always use  Krein spaces. In particular it is crucial to be able to canonically attach a semi-Riemannian spectral triple, hence a Krein space, to a semi-Riemannian spin manifold. However, if the Krein space of spinors depends on the space/time splitting, this is not possible without an extra piece of data.

In this article we  first recall the general definitions concerning hermitian forms on the spinor bundle and Krein spaces in section 2. Then we   give an explicit example of the dependence of the completion on the space/time splitting in section 3. In section 4 we   show that in the case of Lorentzian manifolds, two splittings give rise to equivalent norms exactly when the relative Doppler shift between the two timelike subbundles is bounded on the manifold. In section 5 we extend this result to the other signatures. In the way we are led to define a generalization of the Doppler shift factor as a measure of the relative positions of two subspaces of maximal dimension on which the metric is negative definite. Finally we    sketch two possible ways out to  this non-uniqueness issue in Noncommutative Geometry, leaving their exploration for the future.

\section{Setting}
Let $(M,g)$ be a semi-Riemannian manifold of even dimension and signature $(p,q)$ and suppose that there exists a spinor bundle, i.e. a complex vector bundle $S$ with a left action of the Clifford bundle $Cl(M,g)$.  By a \emph{spinor metric} on $S$, we mean a smooth field of non-degenerate hermitian forms $x\mapsto H_x$ on $S$ such that
\be
H_x(v\cdot\psi,\phi)=H_x(\psi,v\cdot \phi)\label{spinormetric}
\ee
for all $v\in T_xM$, $\psi,\phi\in S_x$.  The existence of a spinor metric was discussed for the first time in \cite{baum}, section 3.3.1 (in German). In \cite{part1}, theorem 3, it is proven that a spinor metric exists iff $(M,g)$ is time-orientable when $p$ and $q$ are even, and iff $(M,g)$ is space orientable when $p$ and $q$ are odd. Given a spinor metric $H$, one can define a non-degenerate product $(.,.)$ on smooth spinor fields with compact support by
\be
(\Psi,\Phi):=\int_MH_x(\Psi(x),\Phi(x)){\rm vol}_g\label{can}
\ee
for all $\Psi,\Phi\in\Gamma_c(S)$. We will call $(.,.)$ the \emph{Krein product}. There exists a multi-vector $v$ ($v\in\Lambda^q TM$ in the even case,  $v\in\Lambda^pTM$ in the odd case) and an integer $r$ such that the operator $\eta$ on $\Gamma_c(S)$ defined by
\be
\eta\Psi(x):=i^rv_x\cdot\Psi(x)\label{fd}
\ee
where $\cdot$ is the natural action of $\Lambda T_xM\simeq Cl(T_xM,g_x)$ on $S_x$, has the following properties:
\begin{enumerate}
\item $\eta^2=1$,
\item $\eta$ is self-adjoint with respect to the Krein-product,
\item $\bra .,.\ket_\eta:=(.,\eta.)$ is positive-definite.
\end{enumerate}
Let us now recall some definitions which can be found in \cite{bognar} (p 49, 52, and 100). Let ${\mathcal K}$ be a complex vector space equipped with a nondegenerate hermitian form $(.,.)$, called the Krein product. A decomposition:
\be
\mathcal{K} = \mathcal{K}_{+} \oplus \mathcal{K}_{-}\label{FundamentalDecomposition}
\ee
is said to be a \emph{fundamental decomposition} of $\mathcal{K}$ if and only if:
\begin{itemize}
	\item The Krein product is positive definite on $\mathcal{K}_{+}$,
	\item The Krein product is negative definite on $\mathcal{K}_{-}$,
	\item The subspaces $\mathcal{K}_{+}$ and $\mathcal{K}_{-}$ are mutually orthogonal with respect to $(\cdot, \cdot)$. 
\end{itemize}
If such a decomposition exists, ${\cal K}$ is said to be \emph{decomposable}. An operator $\eta$ on $\mathcal{K}$ is said to be a \emph{fundamental symmetry} if there exists a fundamental decomposition \eqref{FundamentalDecomposition} such that:
\begin{itemize}
	\item $\eta = 1$ on $\mathcal{K}_{+}$
	\item $\eta = -1$ on $\mathcal{K}_{-}$
\end{itemize}
Clearly fundamental decompositions are in one-to-one correspondence with fundamental symmetries. Note that the restriction of $\pm(.,.)$ on $\mathcal{K}_\pm$ induces a norm. If there exists a fundamental decomposition such that  both $\mathcal{K}_\pm$ are complete with respect to this norm, then we say that $({\cal K},(.,.))$ is a \emph{Krein space}. A crucial result is the equivalence of all fundamental symmetries in a Krein space. More precisely (see \cite{bognar}, cor. 1.2, p. 100):
\begin{theorem} If $({\cal K},(.,.))$ is decomposable and non-degenerate, then it is a Krein space iff, for every fundamental symmetry $\eta$, the scalar product $\bra .,.\ket_\eta$ turns ${\cal K}$ into a Hilbert space (i.e. it is complete). 
\end{theorem}
One can also show that the norms associated to all the fundamental symmetries are equivalent (\cite{nadir}, theorem 2.1). Hence they define a common topology, called by Bognar the \emph{strong topology}. Finally, we recall for future use the following terminology: a subspace of a Krein space is called \emph{maximal negative definite} iff the Krein product is negative definite on it and it is maximal for inclusion among subspaces having this property. 

Let us come back to the case of spinor fields. The operator \eqref{fd} is clearly a fundamental symmetry on $\Gamma_c(S)$. If we define $\Gamma_\eta(S)$ to be the completion of $\Gamma_c(S)$ with respect to $\bra.,.\ket_\eta$, then $(\Gamma_\eta(S),(.,.))$ is a Krein space: this is the content of \cite{baum}, Satz 3.16. However, even though the fundamental symmetries on $\Gamma_\eta(S)$ are all equivalent as explained above, the fundamental symmetries on $\Gamma_c(S)$ are not, and can give rise to different completions. We give an example in the next section.

\section{Counterexample}
Let $\RR^{1,3}$ be Minkowski space, and $S=\RR^{1,3}\times \CC^4$ be the trivial spinor bundle with representation $\rho_x : e_\mu\mapsto \gamma_\mu$, where $e_\mu$ is the canonical basis of $T_x \RR^{1,3}=\RR^{1,3}$ and $\gamma_\mu$ is a fixed set of gamma matrices such that $\gamma_0^\dagger=\gamma_0$, $\gamma_i^\dagger=-\gamma_i$, $i=1,\ldots,3$. For instance we can choose
\be
\gamma_0=\pmatrix{1&0&0&0\cr 0&1&0&0\cr 0&0&-1&0\cr 0&0&0&-1},\gamma_1=\pmatrix{0&0&0&1\cr 0&0&1&0\cr 0&-1&0&0\cr -1&0&0&0}
\ee
as well as $\gamma_2=\pmatrix{0&0&0&-i\cr 0&0&i&0\cr 0&i&0&0\cr -i&0&0&0}$, $\gamma_3=\pmatrix{0&0&1&0\cr 0&0&0&-1\cr -1&0&0&0\cr 0&1&0&0}$. The spinor metric at $x$  is $H_x(\psi,\phi)=\psi^\dagger\gamma_0\phi$. The Krein product on compactly supported section of $S$ is $(\Psi,\Phi)=\int_{\RR^{1,3}}H_x(\Psi(x),\Phi(x))dx^0\ldots dx^3$.

We define the following section of $\End(S)$:
\bea
\eta(x)&:=&\gamma_0\cr
n(x)&:=&\cosh(x_0)\gamma_0+\sinh(x_0)\gamma_1
\eea
Since $\eta(x)$ and $n(x)$ are the images under $\rho_x$ of future-directed timelike vectors, we have that $(\Psi,\eta\Psi)>0$ and $(\Psi,n\Psi)>0$ for all non-zero spinor field with compact support, and $\eta^2=n^2=1$. Let us  notice that
\be
(\Psi,n\Psi)=\bra \Psi,\eta n\Psi\ket_\eta=\int \Psi^\dagger(x) (\cosh x_0+(\sinh x_0)\gamma_0\gamma_1)\Psi(x) dx^0\ldots dx^3\label{eq7}
\ee
Let $\epsilon_1$ be the first vector of the canonical basis of $\CC^4$ and $\phi_k$, $k\in \NN$, be a Dirac delta sequence of  smooth non-negative functions with compact support. Hence $\int_{\RR^4}\phi_k=1$ for all $k$ and $\lim_{k\rightarrow \infty}\int_{\RR^4}f(x)\phi_k(x)dx=f(0)$ for all smooth $f$ with compact support. We can also suppose that $\phi^{1/2}_k$ is smooth, using a Gaussian bump for instance. Then let $y\in \RR^{1,3}$ be some point, and define $\Psi^y_k\in \Gamma_c(S)$ by $\Psi^y_k(x)=\phi_k^{1/2}(x-y)\epsilon_1$. On the one hand we have:
\bea
\|\Psi^y_k\|_\eta^2&=&\int_{\RR^{1,3}} \phi_k(x-y)\epsilon_1^\dagger\gamma_0^2\epsilon_1 dx^0\ldots dx^3\cr
&=&\int_{\RR^{1,3}} \phi_k(x-y)dx^0\ldots dx^3=1, 
\eea
and on the other hand, from \eqref{eq7}:
\bea
\|\Psi^y_k\|_n^2&=&\int_{\RR^{1,3}} \phi_k(x-y))\epsilon_1^\dagger(\cosh x_0+(\sinh x_0)\gamma_0\gamma_1)\epsilon_1 dx^0\ldots dx^3\cr
&=&\int_{\RR^{1,3}} \phi_k(x-y)\cosh x_0 dx^0\ldots dx^3\rightarrow\cosh y_0\mbox{ as }k\rightarrow\infty
\eea
We thus obtain that the two norms are non-equivalent. This does not immediately entail that the completions are different. Indeed, it is well-known that a given Banach space $X$ of infinite dimension can be equipped with two non-equivalent norms for which it is complete. Hence the completion of $X$, which is $X$, is the same for the two norms. However this is a moot point for several reasons. First, finding such an example of a topological vector space which is complete for two non-equivalent norms necessarily involves the axiom of choice (\cite{hba}, cor. 27.47). Moreover, it not hard in the case at hand to prove directly that the completions are different: consider a smooth spinor field $\Psi$ which is of compact support in the spatial directions and decays as $e^{-|x_0|/2}$ as $|x_0|\rightarrow\infty$. An immediate computation shows that $\Psi$ has a finite $\eta$-norm and an infinite $n$-norm. This shows (multiplying with a smooth cut-off function with compact support) that $\Psi$ is the limit of a Cauchy sequence for the $\eta$-norm but not for the $n$-norm. Finally,  and perhaps more importantly, we are not really interested in the question whether the completions are the same set, but rather in knowing whether there exists an isomorphism of Krein spaces from $\Gamma_\eta(S)$ to $\Gamma_n(S)$ covering the identity. In this point of view it is immediate that the identity must be an homeomorphism from $(\Gamma_c(S),\|\ \|_\eta)$ to $(\Gamma_c(S),\|\ \|_n)$, i.e. that the $\eta$-norm and $n$-norm must be equivalent. Thus we see that the equivalence of norms is what really matters. We will focus on this question in the next sections.
 
\section{Condition for the equivalence of spinor norms on a Lorentzian manifold}
Let $n_1=\rho(v_1)$ and $n_2=\rho(v_2)$ be two fundamental symmetries of an (anti-)Lorentzian manifold. We would like to have a constant $k>0$ such that $\|\Psi\|_{n_1}\le k \|\Psi\|_{n_2}$ for all $\Psi\in\Gamma_c(S)$. As above, $\|\Psi\|_{n_1}^2=\bra \Psi,n_2 n_1\Psi\ket_{n_2}$. Now $T=n_2n_1$ is a positive operator for $\bra .,.\ket_{n_2}$ and $(n_2 n_1)^2=n_2(-n_2n_1+2g(v_1,v_2))n_1=-1+2g(v_1,v_2)n_2n_1$. Thus $T$ satisfies the equation 
\be 
T^2-2g(v_1,v_2)T+1=0
\ee
 At any point $x$ of the manifold, the operator $T(x)=n_2n_1(x)$ has therefore eigenvalues $\lambda_\pm=g(v_1,v_2)\pm\sqrt{g(v_1,v_2)^2-1}$. Since $g(v_1,v_2)>0$, the operator norm of $T(x)$ is thus
\be
\|T(x)\|_{n_2}=g(v_1,v_2)+\sqrt{g(v_1,v_2)^2-1}\label{normT}
\ee
We have 
\bea
\|\Psi\|_{n_1}^2&=&\int_M\bra \Psi(x),T(x)\Psi(x)\ket_{n_2}{\rm vol}_g\cr
&\le&\int_M\|T(x)\|_{n_2}\bra \Psi(x),\Psi(x)\ket_{n_2}{\rm vol}_g\cr
&\le&\sup_{M}\|T(x)\|_{n_2}\times\|\Psi\|_{n_2}^2\nonumber
\eea
where the sup is allowed to be infinite. Thus we see that $k$ exists if the expression on the RHS of \eqref{normT} is bounded on the manifold. Conversely, if it is unbounded we can easily build a spinor field with unit $n_2$-norm and arbitrarily large $n_1$-norm by choosing it to be peaked at an eigenspinor of $T(x)$ at a point where $T(x)$ has a very large $n_2$-norm. Since \eqref{normT}  is symmetrical in $v_1,v_2$, we conclude that if $k$ exists then $\|.\|_{n_1}$ and $\|.\|_{n_2}$ are equivalent.

Now, in the 2-plane generated by $(v_1,v_2)$, there is a unique boost which sends $v_1$ to $v_2$. In a pseudo-orthonormal basis $(v_1,v_1^\perp)$ of this 2-plane, this boost has the form 
\be
\pmatrix{\cosh \xi&\sinh \xi\cr \sinh \xi&\cosh \xi}\label{boost}
\ee
where $\xi$ is the rapidity and can be chosen to be positive by the choice of orientation of $v_1^\perp$. The expression \eqref{normT} is then nothing but $e^\xi$, the relative Doppler-shift factor between the two reference frames. We thus obtain

\begin{propo}\label{acc} Let $v_1,v_2$ be two timelike future-directed vector fields on $M$, with $g(v_1,v_1)=g(v_2,v_2)=1$, and define $$\dsf(x):=g(v_1,v_2)+\sqrt{g(v_1,v_2)^2-1}$$
to be the Doppler-shift factor between $v_1$ and $v_2$ at $x$. Then the two fundamental symmetries $n_1=\rho(v_1)$ and $n_2=\rho(v_2)$ define equivalent norms on $\Gamma_c(S)$ iff $k=\sup_{x\in M}\dsf(x)<\infty$. In this case we have $\|\Psi\|_{n_1}\le k^{1/2}\|\Psi\|_{n_2}$ and  $\|\Psi\|_{n_2}\le k^{1/2}\|\Psi\|_{n_1}$ for any $\Psi\in\Gamma_c(S)$. 
\end{propo}

\section{Generalization to all signatures}
\subsection{Doppler shift factor between maximal negative definite subpaces}
In order to generalize proposition \ref{acc} to other signatures, we first need to generalize the Doppler shift factor, or equivalently the hyperbolic angle $\xi$. Let $(V,g)$ be a vector space of dimension $n$ equipped with a nondegenerate bilinear form $g$, and let $V_{i}\oplus V_{i}^\perp$, $i=1,2$, be two fundamental decompositions of the tangent space such that $g$ is negative definite on $V_{i}$, $i=1,2$. We call $s_i$ the $g$-orthogonal symmetry with respect to $V_i$, and $g_i$ the \emph{positive definite} bilinear form defined by $g_i(u,v):=g(s_i(u),v)$, $u,v\in V$. We denote by $\|\ \|_{g_i}$ the operator norm on $\End(V)$ with respect to the scalar product $g_i$.

\begin{lemma}\label{previouslemma} Let  $\Lambda,\Lambda'\in SO(g)$ be such that $\Lambda V_1=\Lambda'V_1=V_2$. Then $\|\Lambda'\|_{g_2}=\|\Lambda\|_{g_2}$.
\end{lemma}
\begin{demo}
Indeed,  $\Lambda'$ can be written $O\Lambda$ with $O\in SO(g)$ such that $O(V_2)=V_2$. As a result, $O\in SO(g_2)$ since it commutes with $s_2$.
\end{demo}

It follows from this lemma that $\|\Lambda\|_{g_2}$ is a measure of the relative position of $V_1$ and $V_2$. In the Lorentzian case  $\Lambda$ can be taken to be the boost \eqref{boost} and we recover the Doppler shift factor. We are hence led to the following definition:

\begin{definition} Let $V_1$ and $V_2$ be two maximal definite negative subspaces. The \emph{Doppler shift factor between $V_1$ and $V_2$} is defined to be  
$$\dsf(V_1,V_2):=\|\Lambda\|_{g_2}$$
where $\Lambda$ is any element of $SO(g)$ such that $\Lambda(V_1)=V_2$.
\end{definition}

Let the ``angle $(V_1,V_2)$'' be the class of the special pseudo-orthogonal transformations which send $V_1$ to $V_2$. We now define a particularly nice representative of this class. Let $\Lambda$ be one of these transformations, and let $*$ stand for the conjugation with respect to $g_2$. We have $A^*=s_2A^\times s_2$ for any operator $A$. In particular $\Lambda^*\Lambda= s_2\Lambda^{-1}s_2\Lambda=s_2s_1$. It follows that $\Lambda^*\Lambda\in SO(g)$. In a Euclidean $2$-plane the product of two orthogonal symmetries is the rotation of twice the angle between the two axes of symmetries. Inspired by this, we define $L$ to be the positive square root of $\Lambda^*\Lambda$ with respect to $g_2$. We thus have  the polar decomposition of $\Lambda$:
\be
\Lambda=OL,\mbox{ with }L^*=L,\ L\ge 0,\  O\in SO(g_2)
\ee
It is a known fact that both $L$ and $O$ belong to $SO(g)$ since $\Lambda$ does (see for instance \cite{MMT}, theorem 5.1). Thus $O\in SO(g)\cap SO(g_2)\simeq SO(V_2)\times SO(V_2^\perp)$, and in particular $O(V_2)=V_2$. It follows that $L(V_1)=V_2$. Let $r(L)$ be the spectral radius of $L$. Since $L^*=L$, we obtain:
\be
\dsf(V_1,V_2)=\|L\|_{g_2}=r(L)\label{srL}
\ee
{\bf Remark:} The angle $(V_2,V_1)$ has $L^{-1}$ as a canonical representative. Using this fact and lemma \ref{lem4} below it follows that the Doppler shift factor is independent of the order in which the pair $(V_1,V_2)$ is written.

\subsection{The operator norm of spin lifts of pseudo-orthogonal transformations}
Let us now suppose that $n$ is even and that $S$ is an irreducible representation space of $Cl(V,g)$. Since we have an isomorphism $\rho : \CC l(V,g)\rightarrow \End(S)$ we identify these two algebras and consider $Spin(g)$ to be a subgroup of $\End(S)$. Whenever we have an element $A$ of $SO(g)$ we will write $\tilde A$ for one of its two lifts to $Spin(g)$. The Krein product on $S$ will be denoted by $H$. The splittings $V=V_i\oplus V_i^\perp$ determine fundamental symmetries $n_i$ on $S$. The associated scalar product, norm and operator norms are written with a $n_i$ index.  Let us recall that $n_i$ can be written $i^rv_1\ldots v_k$, where $r=0$ or $1$, and $v_1,\ldots,v_k$ is a pseudo-orthonormal basis of $V_i$ if $p/q$ is even, of $V_i^\perp$ if $p/q$ is odd. In any case the resulting $n_i$ does not depend on the choice of the pseudo-orthonormal basis. 

\begin{lemma} Let  $\Lambda,\Lambda'\in SO(g)$ be such that $\Lambda V_1=\Lambda'V_1=V_2$. Then $\|\tilde\Lambda'\|_{n_2}=\|\tilde\Lambda\|_{n_2}$.
\end{lemma}
\begin{demo}
Writting $\Lambda'=O\Lambda$ as in lemma \ref{previouslemma},   we have:
\bea
\bra\tilde O\psi,\tilde O\psi\ket_{n_2}&=&H(\tilde O\psi,n_2\tilde O\psi)\cr
&=&H(\psi,\tilde O^{-1}n_2\tilde O\psi),\mbox{ since }\tilde O\in Spin(g)\cr
&=&H(\psi,n_2\psi),\mbox{ since }O(V_2)=V_2\cr
&=&\bra \psi,\psi\ket_{n_2}
\eea
Hence $\tilde O$ is an isometry for the $n_2$-norm, and since $\tilde{\Lambda'}=\pm\tilde O\tilde\Lambda$ we obtain $\|\tilde \Lambda'\|_{n_2}=\|\tilde \Lambda\|_{n_2}$.
\end{demo}

We will need a relation between $\|\tilde\Lambda\|_{n_2}$ and $\|\Lambda\|_{g_2}$. In the Lorentzian case it can be shown that $\|\tilde\Lambda\|_{n_2}=\|\Lambda\|_{g_2}^{1/2}$, which gives back the factor $k^{1/2}$ in proposition \ref{acc}. In the general case   we will only obtain the following estimate:

\begin{propo}\label{propo2}
We have $\|\Lambda\|_{g_2}^{1/2}\le\|\tilde\Lambda\|_{n_2}\le \|\Lambda\|_{g_2}^{\min(p,q)/2}$.
\end{propo}

In order to prove the proposition, we first replace $\Lambda$ with the canonical representative $L$ defined in the previous subsection, so that $\|\Lambda\|_{g_2}=r(L)$. It turns out that  we also have $\|\tilde L\|_{n_2}=r(\tilde L)$. To prove this claim it suffices to show that $\tilde L$ is a normal operator. Let us still denote by $s_2$ the canonical extension of $s_2$ to $\CC l(V,g)$. The anti-linear anti-automorphism $a\mapsto a^*:=s_2(a^T)$ of the Clifford algebra is sent by the spinor representation $\rho$ to the adjunction with respect to the $n_2$-scalar product, i.e. $\rho(a^*)=\rho(a)^*$ (see \cite{part1}, proposition 9). Now let $v\in T_xM$. We have:
\bea
\tilde L^* v\tilde (L^*)^{-1}&=&s_2(\tilde L^T)vs_2(\tilde L^T)^{-1}\cr
&=&s_2(\tilde L^T s_2(v)(\tilde L^T)^{-1})\cr
&=&s_2(\tilde L^{-1} s_2(v)\tilde L)^{T})\cr
&=&s_2(L^{-1}[s_2(v)]),\mbox{ since }w^T=w\mbox{ for }w\in T_xM\cr
&=&s_2s_1L^{-1}(v)\cr
&=&L^2L^{-1}(v)=L(v)\cr
&=&\tilde L v\tilde L^{-1}\nonumber
\eea
Thus $Ad_{\tilde L^*}$ and $Ad_{\tilde L}$ coincide on $V$, and it follows that they coincide on the whole Clifford algebra. Therefore $\tilde L^*$ and $\tilde L$ are proportional, and in particular $\tilde L$ is a normal operator.   It now remains to estimate $r(\tilde L)$. We will need two classical results.
\begin{lemma}\label{lem3} Let $A\in M_N(\CC)$ be an invertible diagonalizable matrix. Then $\mu$ is an eigenvalue of $Ad_A$ iff $\mu$ is the quotient of two eigenvalues of $A$.
\end{lemma}
The proof is immediate by reducing to the case where $A$ is diagonal.
\begin{lemma}\label{lem4} Let $L$ be a pseudo-orthogonal (resp. pseudo-unitary) matrix for a quadratic (resp. sesquilinear) form of signature $(p,q)$. Then the spectrum of $L$ is invariant under $z\mapsto {1\over \bar z}$. Moreover the multiplicities of $\lambda$ and $\bar{\lambda}^{-1}$ as eigenvalues of $L$ are the same.
\end{lemma}
For a proof, see \cite{goh}, Prop. 4.3.3. It follows in particular from this lemma that the number of  eigenvalues of $L$ which have a modulus larger than $1$ cannot exceed $\min(p,q)$, counting multiplicities.
 
Since $L^*=L$, it is diagonalizable in a $g_2$-orthonormal basis $e_1,\ldots,e_n$ of $V$, with eigenvalues $\lambda_1,\ldots,\lambda_n$. The elements $e_I:=e_{i_1}\ldots e_{i_k}$ where $\{i_1;\ldots;i_k\}$ is a subset of $\{1;\ldots;n\}$ written in increasing order form a basis of the Clifford algebra, and one has $Ad_{\tilde L}(e_I)=Le_{i_1}\ldots Le_{i_k}=\lambda_{i_1}\ldots\lambda_{i_k}e_I$. Thus the spectrum of $Ad_{\tilde L}$ is:
\be
\sigma(Ad_{\tilde L})=\{\lambda_{i_1}\ldots\lambda_{i_k}|1\le i_1<\ldots<i_k\le n\}\label{spectre1}
\ee
and by lemma \ref{lem3}, we also have
\be
\sigma(Ad_{\tilde L})=\{{\tilde \lambda\over \tilde\lambda'}|\tilde\lambda,\tilde\lambda'\in \sigma(\tilde L)\}\label{spectre2}
\ee
Now   $\tilde L$ is pseudo-unitary, hence by  lemma \ref{lem4} and \eqref{spectre2} we obtain
\be
r(Ad_{\tilde L})=r(\tilde L)^2\label{eqtruc}
\ee
Now from \eqref{spectre1} we see that $r(Ad_{\tilde L})$ is the maximal value of $|\lambda_{i_1}\ldots\lambda_{i_k}|$, which is reached when all  the eigenvalues  of $L$ which are larger than $1$ are multiplied together. Hence we obtain
\be 
r(Ad_{\tilde L})\le r(L)^{\min(p,q)}\label{eqmachin}
\ee
Putting \eqref{eqtruc} and \eqref{eqmachin} together we obtain 
\be
\|\tilde\Lambda\|_{n_2}=\|\tilde L\|_{n_2}=r(\tilde L)\le r(L)^{\min(p,q)\over 2}=\|L\|_{g_2}^{\min(p,q)\over 2}=\|\Lambda\|_{g_2}^{\min(p,q)\over 2}
\ee
Moreover, we clearly have $r(Ad_{\tilde L})\ge r(L)$ hence $r(\tilde L)\ge r(L)^{1/2}$ from which we obtain the other half of the sought after estimate. Some remarks are in order.

\begin{itemize}
\item In Lorentzian or anti-Lorentzian signature, the estimate can be replaced be the equality $\|\tilde \Lambda\|_{n_2}=\|\Lambda\|_{g_2}^{1/2}$.
\item In Euclidean or anti-Euclidean signature, the splitting is trivial and the estimate degerates into $\|\tilde \Lambda\|=\|\Lambda\|^{1/2}=1$, which is obvious. 
\item In other signatures the estimate is optimal. Indeed, let the metric be $g={\rm diag}(1,-1,\ldots,1,-1,1,\ldots,1)$, where the couple $(1,-1)$ appears $q$ times. Let  $V_2$ be such that $g_2$ is the canonical scalar product on $\RR^n$ and $\Lambda$ be the ``$q$-boost'' $\Lambda=\pmatrix{\cosh x_1&\sinh x_1\cr \sinh x_1&\cosh x_1}\oplus\ldots\oplus \pmatrix{\cosh x_q&\sinh x_q\cr \sinh x_q&\cosh x_q}\oplus I_{p-q}$, with $0<x_1<\ldots<x_q$. Then $L=\Lambda$ and $r(Ad_{\tilde L})=\exp(\sum_{k=1}^qx_k)$, thus $\|\tilde \Lambda\|_{n_2}=\exp({\sum_{k=1}^q x_k\over 2})$ while $\|\Lambda\|_{g_2}=\exp x_q$. The estimate thus yields $\exp{x_q\over 2}\le \exp{\sum_{k=1}^q x_k\over 2}\le \exp{qx_q\over 2}$. We see that  both side can degenerate into an equality. 
\end{itemize}

\subsection{Condition for the equivalence of spinor norms}
We now consider a spinor bundle $S$ over a spin semi-Riemannian manifold $(M,g)$ of even dimension. We consider $n_1=\rho(i^r v_1)$ and $n_2=\rho(i^rv_2)$ two fundamental symmetries of $\Gamma_c(S)$, where this time $v_1$ and $v_2$ are multivector fields. At each point $x\in M$, $v_1(x)$ and $v_2(x)$ generate subspaces of $T_xM$ of definite signatures. If $p$ or $q$ is even, the subspaces generated by the $v_i(x)$ are negative, and we call them $V_i(x)$. If $p$ or $q$ is odd, they are positive and we call them $V_i(x)^\perp$.  Hence $V_{i}(x)\oplus V_{i}(x)^\perp$, $i=1,2$, are two fundamental decompositions of the tangent space as in the previous subsection. For each $x$ we let $\Lambda(x)$ denote any element of $SO(T_xM)$ such that $\Lambda(x)V_1(x)=V_2(x)$. From this we obtain $\tilde\Lambda (x)n_1(x)\tilde \Lambda(x)^{-1}=n_2(x)$, and it follows that for any spinor $\psi\in S_x$:
\bea
\|\psi\|_{n_1}^2&=&H(\psi,n_1\psi)\cr
&=&H(\psi,\tilde\Lambda^{-1} n_2\tilde \Lambda\psi)\cr
&=&\|\tilde\Lambda\psi\|_{n_2}^2\cr
\Rightarrow \|\psi\|_{n_1}&\le&\|\tilde \Lambda\|_{n_2}\|\psi\|_{n_2}\label{estim}
\eea
where $n_1,n_2,\tilde\Lambda$ and $H$ are evaluated at $x$. Let us remark that even if a smooth field $x\mapsto \Lambda(x)$ is not defined globally on $M$, the map $x\mapsto\|\tilde\Lambda(x)\|_{n_2}$ is, since at every $x$ this number only depends on the angle $(V_1(x),V_2(x))$, as  defined  in previous subsection.  Integrating on $M$, and using $\|\tilde\Lambda^{-1}\|_{n_1}=\|\tilde\Lambda\|_{n_2}$, we  obtain that   if $x\mapsto\|\tilde\Lambda(x)\|_{n_2}$ is bounded then the $n_1$-norm and the $n_2$-norm are equivalent. It is easy to prove the converse. Indeed, let $\psi\in S_{x_0}$ be such that $\|\tilde\Lambda(x_0) \psi\|_{n_2}=\|\tilde\Lambda(x_0)\|_{n_2}\|\psi\|_{n_2}$. Let $U$ be a small open set around $x_0$ such that the spinor bundle is trivial over $U$. Thus $S_{|U}\simeq U\times S_{x_0}$. We let $f_k$ be a sequence of smooth functions with support in $U$ converging to the delta function $\delta_{x_0}$, and define the spinor field $\Psi_k$ to be $\Psi_k(x)=f_k(x)\psi$ on $U$ and $0$ outside $U$. Then the quotient ${\|\Psi_k\|_{n_1}\over \|\Psi_k\|_{n_2}}$ is easily seen to converge to $\|\tilde\Lambda(x_0)\|_{n_2}$, which is as large as we want by hypothesis. We thus have obtained that the $n_1$ and $n_2$-norms are equivalent iff $x\mapsto \|\tilde \Lambda(x)\|$ is bounded on $M$. Hence proposition \ref{propo2} yields the following theorem:

\begin{theorem} Let $(M,g)$ be a semi-Riemannian spin manifold of even dimension and $TM=V_1\oplus V_1^\perp=V_2\oplus V_2^\perp$ be two splittings of the tangent bundle into negative and positive definite subbundles respectively. Let $n_1$ and $n_2$ be the associated fundamental symmetries of the indefinite inner product space $\Gamma_c(S)$.  Then the $n_1$-norm and the $n_2$-norm are equivalent iff $\sup_{x\in M}\dsf(V_1(x),V_2(x))<\infty$.
\end{theorem}

\section{Doppler classes and the consequences for Noncommutative Geometry}
As we have said in the introduction, in Semi-Riemannian Noncommutative Geometry one would like to define a canonical ``spectral triple'' out of the data of a spin manifold. The above arguments show that this cannot be done: one must provide an additional piece of information. In order to be more specific, let us consider the case of a Lorentzian spin manifold $(M,g)$. We will say that two future-directed timelike vector fields $v_1,v_2$ are in the same \emph{Doppler class} if $g(v_1,v_2)$ is bounded on $M$. By proposition \ref{acc}, the extra piece of information which is needed in order to define a Krein space of spinor fields on $M$ is a Doppler class. A possibility that would seem to remain open would be to single out a particularly natural Doppler class on a given manifold. Let us show by some examples that one should not cling too much to that hope. The most natural choice would be the class of timelike geodesic vector fields, could it be proven that they all belonged to the same Doppler class. This is indeed true in Minkowski spacetime of dimension $1+2$, but unfortunately this is false in dimension $1+3$, as the following example shows:
\bea
v_0&=&\cosh x_3\cr
v_1&=&\sinh x_3\cr
v_2=v_3&=&0
\eea
with coordinates $(x_0,\ldots,x_3)$, defines a normalized future-directed timelike vector field $v$ the  integral curves  of which  are  straight lines lying in $2$-planes of  constant $x_2,x_3$. Moreover   $g(v,e_0)$ is undounded, with $e_0=(1,0,0,0)$. Hence $v$ and $e_0$ do not belong to the same Doppler class.  Setting gravity in does not seem to help: the radially ingoing and outgoing geodesic vector fields in the Schwarzschild spacetime can be proven not to be in the same Doppler class. In order to avoid these pathologies, we have to put more restrictions on the vector fields.  Covariantly constant timelike vector fields would fit the bill: if $v_1,v_2$ are two such fields, then $\partial_\mu g(v_1,v_2)=g(\nabla_\mu v_1,v_2)+g( v_1,\nabla_\mu v_2)=0$, thus $v_1$ and $v_2$ belong to the same Doppler class. It is interesting to note that such vector fields seem to play an important role in discrete spectral spacetimes (\cite{SST2}). However these vector fields seldom exist. For instance, it can be shown that the only four dimensional vacuum solution of the Einstein equation on which non-null covariantly constant vector fields exist is the flat solution \cite{batista}.

We conclude from these examples that the best strategy for Noncommutative Semi-Riemannian Geometry is probably \emph{not to} use Krein spaces, but ``pre-Krein spaces'', that is non-degenerate and decomposable inner product spaces, the canonical example of which being the  space $\Gamma_c(S)$  equipped with the  form  \eqref{can}. All its possible fundamental symmetries and completions would then have to be equally considered.

\bibliographystyle{unsrt}
\bibliography{../generalbib/SSTbiblio}

\end{document}